\documentclass{amsart}

\input xy
\xyoption{all}

\newcommand{\w}{\omega}

\newcommand{\id}{\mathrm{id}}
\newcommand{\Ra}{\Rightarrow}
\newcommand{\U}{\mathcal U}
\newcommand{\W}{\mathcal W}
\newcommand{\cl}{\operatorname{cl}}
\newcommand{\pr}{\operatorname{pr}}
\newcommand{\IN}{\mathbb{N}}
\newcommand{\diam}{\operatorname{diam}}
\newcommand{\Rea}{\mathbb{R}}

\newtheorem{theorem}{Theorem}[section]
\newtheorem{proposition}[theorem]{Proposition}
\newtheorem{corollary}[theorem]{Corollary}
\newtheorem{example}[theorem]{Example}
\newtheorem{lemma}[theorem]{Lemma}
\newtheorem{claim}[theorem]{Claim}
\newtheorem{remark}[theorem]{Remark}
\theoremstyle{definition}
\newtheorem{definition}[theorem]{Definition}

\title{A spectral characterization of skeletal maps}
\author{Taras Banakh, Andrzej Kucharski, and Marta Martynenko, }
\address{T.~Banakh: Faculty of Mechanics and Mathematics, Ivan Franko National University of Lviv (Ukraine) and \\
Instytut Matematyki, Jan Kochanowski University, Kielce (Poland)}
\email{t.o.banakh@gmail.com}
\address{A.~Kucharski: Institute of Mathematics, University of Silesia, ul. Bankowa 14, 40-007 Katowice (Poland)}
\email{akuchar@math.us.edu.pl}
\address{M.~Martynenko: Faculty of Mechanics and Mathematics, Ivan Franko National University of Lviv (Ukraine)}
\email{martamartynenko@ukr.net}

\begin{document}
\begin{abstract}
We prove that a map between two realcompact spaces is skeletal if and only if it is homeomorphic to the limit map of a skeletal morphism  between $\w$-spectra with surjective limit projections.
\end{abstract}
\maketitle

In this paper we present two characterizations of skeletal maps between realcompact topological spaces.
All maps considered in this paper are continuous and all spaces are Tychonoff. For a subset $A$ of a topological space $X$ by $\cl A$ we shall denote the closure of $A$ in $X$.

A map $f:X\to Y$ is called {\em skeletal} if for each nowhere dense subset $A\subset Y$ the preimage $f^{-1}(A)$ is nowhere dense in $X$. This is equivalent to saying that for each non-empty open set $U\subset X$ the closure $\overline{f(U)}$ has non-empty interior in $Y$, see \cite{Miod}.

The latter definition can be localized as follows. A map $f:X\to Y$ between two topological spaces is called
\begin{itemize}
\item {\em skeletal at a point $x\in X$} if for each neighborhood $U\subset X$ of $x$ the closure $\cl_Y(f(U))$ of $f(U)$ has non-empty interior in $Y$;
\item {\em skeletal at a subset} $A\subset X$ if $f$ is skeletal at each point $x\in A$.
\end{itemize}

It is clear that a map $f:X\to Y$ is skeletal if and only if $f$ is skeletal at each point $x\in X$.

\section{Characterizing skeletal maps between metrizable Baire spaces}

It is clear that each open map is skeletal. For closed maps between metrizable Baire spaces this implication can be partly reversed.
Let us recall that a topological space $X$ is {\em Baire} if for any sequence $(U_n)_{n\in\w}$ of open dense subsets $U_n\subset X$ the intersection $\bigcap_{n\in\w}U_n$ is dense in $X$.

We shall say that a map $f:X\to Y$ between topological spaces is
\begin{itemize}
\item {\em open at a point $x\in X$} if for each neighborhood $U\subset X$ of $x$ the image $f(U)$ is a neighborhood of $f(x)$;
\item {\em open at a subset $A\subset X$} if $f$ is open at each point $x\in A$;
\item {\em densely open} if $f$ is open at some dense subset $A\subset X$.
\end{itemize}

It is easy to see that each densely open map is skeletal. The converse is true for skeletal maps between metrizable compacta, and more generally, for closed skeletal maps defined on metrizable Baire spaces.

\begin{theorem}\label{l1} For a closed map $f:X\to Y$ defined on a metrizable Baire space $X$ the following conditions are equivalent:
\begin{enumerate}
\item $f$ is skeletal;
\item $f$ is skeletal at a dense subset of $X$;
\item $f$ is densely open;
\item $f$ is open at a dense $G_\delta$-subset of $X$.
\end{enumerate}
\end{theorem}

\begin{proof} The implications $(4)\Ra(3)\Ra(2)\Ra(1)$ are trivial and hold without any conditions on $f$.

To prove the implication $(1)\Ra(4)$, fix a metric $d$ generating the topology of the metrizable space $X$. For every $n\in\IN$ consider the family $\U_n$ of all non-empty open subsets $U\subset X$ such that $\diam  (U)<1/n$ and $f(U)$ is open in $f$. The skeletal property of $f$ implies that the union $\bigcup\U_n$ is dense in $X$. Since the space $X$ is Baire, the intersection $A=\bigcap_{n=1}^\infty \bigcup\U_n$ is a dense $G_\delta$-set in $X$. It is clear that $f$ is open at the set $A$.
\end{proof}

The following simple example shows that the metrizability of $X$ is essential in Theorem~\ref{l1} and cannot be weakened to the first countability.

\begin{example} The projection $\pr:A\to[0,1]$ from the Aleksandrov ``two arrows'' space $A=[0,1)\times\{0\}\cup(0,1]\times\{1\}$ onto the interval is skeletal. Yet it is open at no point $x\in A$.
\end{example}

\section{Skeletal and densely open squares}

In this section the notions of skeletal and densely open maps are generalized to square diagrams. These generalized properties will be used in the spectral characterization of skeletal maps given in the next section.

\begin{definition} Let $\mathcal D$ be a commutative diagram
$$\xymatrix{
\tilde X\ar[r]^{\tilde f}\ar[d]_{p_X}&\tilde Y\ar[d]^{p_Y}\\
X\ar[r]_{f}&Y
}$$ consisting of continuous maps between topological spaces. The commutative square $\mathcal D$ is called
\begin{itemize}
\item {\em open at a point $x\in X$} if for each neighborhood $U\subset X$ of $x$ the point $f(x)$ has a neighborhood $V\subset Y$ such that $V\subset f(U)$ and $p^{-1}_Y(V)\subset\tilde f(p_X^{-1}(U))$;
\item {\em open at a subset $A\subset X$} if $\mathcal D$ is open at each point $x\in A$;
\item {\em densely open} if it is open at some dense subset $A\subset X$;
\smallskip

\item {\em skeletal at a point $x\in X$} if for each neighborhood $U\subset X$ of $x$ there is a non-empty open set $V\subset Y$ such that $V\subset \cl f(U)$ and $p_Y^{-1}(V)\subset \cl \tilde f(p_X^{-1}(U))$;
\item {\em skeletal} at a subset $A\subset X$ if $\mathcal D$ is skeletal at each point $x\in A$;
\item {\em skeletal} if $\mathcal D$ is skeletal at $X$.
\end{itemize}
\end{definition}

\begin{remark}\label{rem1} If the square $\mathcal D$ is skeletal (at a point $x\in X$), then the map $f$ is skeletal (at the point $x$).
\end{remark}

\begin{remark} A map $f:X\to Y$ is skeletal (resp. open) at a subset $A\subset X$ if and only if the square  $$\xymatrix{
X\ar[r]^{f}\ar[d]_{\id_X}&Y\ar[d]^{\id_Y}\\
X\ar[r]_{f}&Y
}$$is skeletal (resp. open) at the subset $A$.
\end{remark}

It is easy to see that each densely open square is skeletal. Under some conditions the converse is also true. The following proposition is a ``square'' counterpart of the characterization Theorem~\ref{l1a}.

\begin{proposition}\label{l1a} Let $\mathcal D$ be a commutative diagram
$$\xymatrix{
\tilde X\ar[r]^{\tilde f}\ar[d]_{p_X}&\tilde Y\ar[d]^{p_Y}\\
X\ar[r]_{f}&Y
}$$ consisting of continuous maps between topological spaces such that the map $\tilde f:\tilde X\to\tilde Y$ is closed, the projection $p_Y$ is surjective, and the space $X$ is metrizable and Baire. Then the following conditions are equivalent:
\begin{enumerate}
\item the square $\mathcal D$ is skeletal;
\item $\mathcal D$ is skeletal at a dense subset of $X$;
\item $\mathcal D$ is densely open;
\item $\mathcal D$ is open at a dense $G_\delta$-subset of $X$.
\end{enumerate}
\end{proposition}

\begin{proof} The implications $(4)\Ra(3)\Ra(2)\Ra(1)$ are trivial and hold without any conditions on $\mathcal D$.

To prove the implication $(1)\Ra(4)$, assume that the square $\mathcal D$ is skeletal. First let us prove two auxiliary claims.

\begin{claim}\label{cl1} For each non-empty open  subset $U\subset X$ there is a non-empty open set $V\subset Y$ such that $V\subset f(U)$ and $p^{-1}_Y(V)\subset \tilde f(p^{-1}_X(U))$.
\end{claim}

\begin{proof} Using the regularity of the space $X$, find a non-empty open subset $W\subset X$ whose closure $\overline{W}$ lies in the open set $U$. Since the square $\mathcal D$ is skeletal, for the set $W$ there is a non-empty open set $V\subset Y$ such that $p^{-1}_Y(V)\subset\cl \tilde f(p_X^{-1}(W))$. Taking into account that the map $\tilde f$ is closed, we see that the set $\tilde f(p^{-1}_X(\overline{W}))$ is closed in $\tilde Y$ and hence
$$p^{-1}_Y(V)\subset\cl \tilde f(p_Y^{-1}(W))\subset \tilde f(p^{-1}_X(\overline{W}))\subset \tilde f(p_X^{-1}(U)).$$
Applying to these inclusions the surjective map $p_Y$, we get the inclusion
$$V=p_Y(p_Y^{-1}(V))\subset p_Y\circ\tilde f(p_X^{-1}(U))=f\circ p_X(p^{-1}_X(U))\subset f(U).$$
\end{proof}

\begin{claim}\label{cl2} Each non-empty open set $U\subset X$ contains a non-empty open set $W\subset U$ such that $f(W)$ is open in $Y$ and $\tilde f(p^{-1}_X(W))=p^{-1}_Y(f(W))$.
\end{claim}

\begin{proof} By Claim~\ref{cl1}, there is a non-empty open set $V\subset Y$ such that $V\subset f(U)$
and $p^{-1}_Y(V)\subset \tilde f(p^{-1}_X(U))$. Then the open set $W=U\cap f^{-1}(V)$ has the required properties. Indeed, its image $f(W)=V$ is open in $Y$.
Also the inclusion $p_Y^{-1}(V)\subseteq \tilde f(p_X^{-1}(U))$ implies
$$\begin{aligned} \tilde f(p_X^{-1}(W))&=\tilde f(p_X^{-1}(U\cap f^{-1}(V)))= \tilde f(p_X^{-1}(U)\cap p_X^{-1}(f^{-1}(V)))=\\
&=\tilde f(p_X^{-1}(U)\cap \tilde f^{-1}(p_Y^{-1}(V)))=
\tilde f(p_X^{-1}(U))\cap p_Y^{-1}(V)=p_Y^{-1}(V)=p_Y^{-1}(f(W)).
\end{aligned}
$$\end{proof}

Let $\W$ be the family of all non-empty open sets $W\subset X$ such that $f(W)$ is open in $Y$ and $p^{-1}_Y(f(W))=\tilde f(p^{-1}_X(W))$. Fix any metric $d$ generating the topology of the metrizable space $X$ and for every $n\in\w$ consider the subfamily $\W_n=\{W \in\W:\diam(W)<2^{-n}\}$. By Claim~\ref{cl2}, the union $\bigcup\W_n$ is an open dense subset of $X$. Since $X$ is a Baire space, the intersection $A=\bigcap_{n\in\w}\bigcup\W_n$  is a dense $G_\delta$-set in $X$. To finish the proof, observe that the diagram $\mathcal D$ is open at the dense $G_\delta$-set $A$.
\end{proof}

\section{Skeletal squares and inverse spectra}

In this section we detect morphisms between inverse spectra, inducing skeletal maps between their limit spaces. At first we need to recall some standard information about inverse spectra, see \cite[\S 2.5]{En} and \cite[Ch.1]{Chi} for more details.

For an inverse spectrum $\mathcal S=\{X_\alpha,p_\alpha^\beta,A\}$ consisting of topological spaces and continuous bonding maps, by $$\lim \mathcal S=\{(x_\alpha)_{\alpha\in A}\in\prod_{\alpha\in A}X_\alpha:\forall \alpha\le \beta \;\; p^\beta_\alpha(x_\beta)=x_\alpha\}$$we denote the limit of $\mathcal S$ and by $p_\alpha:\lim \mathcal S\to X_\alpha$, $p_\alpha:x\mapsto x_\alpha$, the limit projections.

Let $\mathcal S_X=\{X_\alpha,p^\beta_\alpha,A\}$ and $\mathcal S_Y=\{Y_\alpha,\pi^\beta_\alpha,A\}$ be  two inverse spectra indexed by the same directed partially ordered set $A$. A {\em morphism} $\{f_\alpha\}_{\alpha\in A}:\mathcal S_X\to\mathcal S_Y$ between these spectra is a family
of maps $\{f_\alpha:X_\alpha\to Y_\alpha\}_{\alpha\in A}$ such that $f_\alpha\circ p^\beta_\alpha=\pi^\beta_\alpha\circ f_\beta$ for any elements $\alpha\le\beta$ in $A$.

Each morphism $\{f_\alpha\}_{\alpha\in A}:\mathcal S_X\to\mathcal S_Y$ of inverse spectra induces a limit map
$$\lim f_\alpha:\lim\mathcal S_X\to\lim \mathcal S_Y,\;\;\lim f_\alpha:(x_\alpha)_{\alpha\in A}\mapsto (f_\alpha(x_\alpha))_{\alpha\in A}$$ between the limits of these inverse spectra.

For indices $\alpha\le\beta$ in $A$ the commutative squares
$$\xymatrix{
\lim\mathcal S_X\ar_{p_\alpha}[d]\ar^{\lim f_\alpha}[r]&\lim\mathcal S_Y\ar^{\pi_\alpha}[d]\\
X_\alpha\ar_{f_\alpha}[r]&Y_\alpha
}
\xymatrix{&\ar@{}[d]|{\mbox{and}}&\\
&&}
\xymatrix{
X_\beta\ar_{p^\beta_\alpha}[d]\ar^{f_\beta}[r]&Y_\beta\ar^{\pi_\alpha^\beta}[d]\\
X_\alpha\ar_{f_\alpha}[r]&Y_\alpha
}
$$
are called respectively the {\em limit $\downarrow_\alpha$-square} and the {\em bonding ${\downarrow}^\beta_\alpha$-square} of the morphism $\{f_\alpha\}$.
\smallskip

We shall say that the morphism $\{f_\alpha\}_{\alpha\in A}:\mathcal S_X\to\mathcal S_Y$ has
\begin{itemize}
\item is {\em skeletal} if each map  $f_\alpha:X_\alpha\to Y_\alpha$, $\alpha\in A$, is skeletal;
\item has {\em skeletal limit squares} if for every index $\alpha\in A$ the limit $\downarrow_\alpha$-square is skeletal;
\item has {\em skeletal bonding squares} if for every indices $\alpha\le\beta$ in $A$ the bonding $\downarrow_\alpha^\beta$-square is skeletal.
\end{itemize}

Our aim is to find conditions on a morphism $\{f_\alpha\}:\mathcal S_X\to\mathcal S_Y$ of spectra implying the skeletality of the limit map $f=\lim f_\alpha:\lim\mathcal S_X\to\lim\mathcal S_Y$.

\begin{proposition}\label{p1} For a morphism $\{f_\alpha\}_{\alpha\in A}:\mathcal S_X\to\mathcal S_Y$ between inverse spectra $\mathcal S_X=\{X_\alpha,p_\alpha^\beta,A\}$ and $\mathcal S_Y=\{X_\beta,\pi_\alpha^\beta,A\}$ with surjective limit projections, the limit map $\lim f_\alpha:\lim \mathcal S_X\to\lim\mathcal S_Y$ is skeletal if the morphism $\{f_\alpha\}$ has skeletal limit squares. \end{proposition}

\begin{proof} We need to show that the limit map $f=\lim f_\alpha:X\to Y$ is skeletal, where $X=\lim \mathcal S_X$, $Y=\lim \mathcal S_Y$. Given any non-empty open set $U\subset X$, we need to find a non-empty open set $V\subset Y$ such that $V\subset\cl f(U)$. By the definition of the topology of the limit space $X=\lim\mathcal S_X$, there is an index $\alpha\in A$ and a non-empty open set $U_\alpha\subset X_\alpha$ such that $U\supset p_{\alpha}^{-1}(U_\alpha)$. Since the limit $\downarrow_\alpha$-square
$$\xymatrix{
X\ar[r]^{f}\ar[d]_{p_\alpha}&Y\ar[d]^{\pi_\alpha}\\
X_\alpha\ar[r]_{f_\alpha}&Y_\alpha}
$$
is skeletal, for the open set $U_\alpha\subset X_\alpha$ there exists a non-empty open set $V_\alpha\subset Y_\alpha$ such that the open set $V=\pi_{\alpha}^{-1}(V_\alpha)$ lies in the closure of the set $f(p^{-1}_\alpha(U_\alpha))$, which lies in the closure of the set $f(U)$.
\end{proof}

It turns out that in some cases the skeletality of squares is preserved by limits.

A partially ordered set $A$ is called {\em $\kappa$-directed} for a cardinal number $\kappa$ if each subset $K\subset A$ of cardinality $|C|\le\kappa$ has an upper bound in $A$.

For a topological space $X$ by $\pi w(X)$ we denote the {\em $\pi$-weight} of $X$, that is, the smallest cardinality $|\mathcal B|$ of a $\pi$-base  $\mathcal B$ for $X$. We recall that a family $\mathcal B$ of non-empty open subsets of $X$ is called a {\em $\pi$-base} for $X$ if each non-empty open subset of $X$ contains a set $U\in\mathcal B$.

\begin{proposition}\label{p1a} Let $\{f_\alpha\}_{\alpha\in A}:\mathcal S_X\to\mathcal S_Y$ be a morphism between inverse spectra $\mathcal S_X=\{X_\alpha,p_\alpha^\beta,A\}$ and $\mathcal S_Y=\{X_\beta,\pi_\alpha^\beta,A\}$ with surjective limit projections. If for some $\alpha\in A$ and the cardinal $\kappa=\pi w(Y_\alpha)$ the index set $A$ is $\kappa$-directed, then the limit $\downarrow_\alpha$-square is skeletal provided that for any $\beta\ge\alpha$ in $A$ the bonding $\downarrow_\alpha^\beta$-square is skeletal.
 \end{proposition}

\begin{proof} Assuming that the limit $\downarrow_\alpha$-square is not skeletal, we can find a non-empty open set $U_\alpha\subset X_\alpha$ such that for any non-empty open set $V_\alpha\subset Y_\alpha$ we get $\pi_\alpha^{-1}(V_\alpha)\not\subset\cl f(U)$ where $U=p_\alpha^{-1}(U_\alpha)$ and $f=\lim f_\alpha$ is the limit map. Fix a $\pi$-base $\mathcal B$ for the space $Y_\alpha$ having cardinality $|\mathcal B|=\pi w(Y_\alpha)\le\kappa$.
For every set $V\in\mathcal B$ the open set $\pi_\alpha^{-1}(V)\setminus\cl f(U)$ is not empty and hence contains a  set of the form $\pi_{\alpha_V}^{-1}(W_V)$ for some index $\alpha_V\ge\alpha$ in $A$ and some non-empty open set $W_V\subset Y_{\alpha_V}$. Since the index set $A$ is $\kappa$-directed, the set $\{\alpha_V:V\in\mathcal B\}$ has an upper bound $\beta\in A$.

By our hypothesis, the bonding $\downarrow_\alpha^\beta$-square
is skeletal. Then for the open subset $U_\beta=(p_\alpha^\beta)^{-1}(U_\alpha)$ of $X_\beta$ we can find a non-empty open set $V\subset Y_\alpha$ such that $(\pi_\alpha^\beta)^{-1}(V)\subset \cl f_\beta(U_\beta)$. We lose no generality assuming that $V\in\mathcal B$. In this case the choice of the set $W_V$ guarantees that $\pi_{\alpha_V}^{-1}(W_V)\subset \pi_{\alpha}^{-1}(V)\setminus f(U)$. Then the open subset $W_\beta=(\pi^{\beta}_{\alpha_V})^{-1}(W_V)=\pi_\beta(\pi_{\alpha_V}^{-1}(W_V))$ of
$ (\pi_\alpha^\beta)^{-1}(V)$
does not intersect the set $\pi_\beta\circ f(U)=f_\beta\circ p_\beta(U)=f_\beta(U_\beta)$ and hence cannot lie in $\cl f_\beta(U_\beta)$.
This contradiction shows that the limit $\downarrow_\alpha$-square is skeletal.
\end{proof}

\begin{corollary}\label{c1a} Let $\{f_\alpha\}_{\alpha\in A}:\mathcal S_X\to\mathcal S_Y$ be a morphism between inverse spectra $\mathcal S_X=\{X_\alpha,p_\alpha^\beta,A\}$ and $\mathcal S_Y=\{X_\beta,\pi_\alpha^\beta,A\}$ with surjective limit projections. If for the cardinal $\kappa=\sup\{\pi w(Y_\alpha):\alpha\in A\}$ the index set $A$ is $\kappa$-directed, then the morphism $\{f_\alpha\}_{\alpha\in A}$ has skeletal limit squares provided it has skeletal bonding squares.
 \end{corollary}

For $\pi\tau$-spectra, Proposition~\ref{p1} can be partly reversed. First let us introduce the necessary definitions.

Let $\tau$ be an infinite cardinal number. We shall say that an inverse spectrum $\mathcal S=\{X_\alpha,p_\alpha^\beta,A\}$ is a {\em $\pi\tau$-spectrum} (resp. a {\em $\tau$-spectrum}) if
\begin{itemize}
\item each space $X_\alpha$, $\alpha\in A$, has $\pi$-weight $\pi w(X_\alpha)\le\tau$ (resp. weight $w(X_\alpha)\le\tau$);
\item the index set $A$ is {\em $\tau$-directed} in the sense that each subset $B\subset A$ of cardinality $|B|\le\tau$ has an upper bound in $A$;
\item the index set $A$ is  {\em $\w$-complete} in the sense that each countable chain $C\subset A$ has the smallest upper bound $\sup C$ in $A$;
\item the spectrum $\mathcal S$ is {\em $\tau$-continuous} in the sense that for any directed subset $C\subset A$ with $\gamma=\sup C$ the limit map $\lim p^\gamma_\alpha:X_\gamma\to\lim \{X_\alpha,p^\beta_\alpha,C\}$ is a homeomorphism.
\end{itemize}

A subset $C$ of a directed poset $A$ is called
\begin{itemize}
\item {\em cofinal} if for any $\alpha\in A$ there is an index $\beta\in C$ with $\alpha\le\beta$;
\item {\em $\tau$-closed} if for each directed subset $D\subset C$ that has the lowest upper bound $\sup D$ in $A$ we get $\sup D\in C$;
\item {\em $\tau$-stationary} if $C$ has non-empty intersection with any cofinal $\tau$-closed subset of $A$.
\end{itemize}

\begin{theorem}\label{p2} Let $\{f_\alpha\}_{\alpha\in A}:\mathcal S_X\to \mathcal S_Y$ be a morphism between two $\pi\tau$-spectra $\mathcal S_X=\{X_\alpha,p^\beta_\alpha,A\}$ and $\mathcal S_Y=\{Y_\alpha,\pi^\beta_\alpha,A\}$ with surjective limit projections. If the limit map $\lim f_\alpha:\lim \mathcal S_X\to\lim\mathcal S_Y$ is skeletal, then for some cofinal $\tau$-closed subset $B\subset A$ the morphism $\{f_\alpha\}_{\alpha\in B}$ is skeletal and has skeletal bonding and limit squares.
\end{theorem}

\begin{proof} To simplify denotations, let $X=\lim \mathcal S_X$, $Y=\lim \mathcal S_Y$, and $f=\lim f_\alpha:X\to Y$. First we show that the set
$$B=\{\alpha\in A:\mbox{the limit $\downarrow_\alpha$-square is skeletal}\}$$is cofinal and $\tau$-closed in $A$. For this we shall prove an auxiliary statement:

\begin{claim}\label{cl:1nn} For every $\alpha\in A$ there is $\beta\in A$, $\beta\ge \alpha$, such that for any non-empty open set
$U\subset X_\alpha$ there is a non-empty open set $V\subset Y_\beta$ such that  $\pi_\beta^{-1}(V)\subseteq\cl f(p_\alpha^{-1}(U)).$
\end{claim}

\begin{proof} In the space $X_\alpha$ fix a $\pi$-base  $\mathcal B$ of cardinality $|\mathcal B|=\pi w(X_\alpha)\le\tau$. For every set $U\in\mathcal B$ the preimage $p_\alpha^{-1}(U)$ is a non-empty open set in $X=\lim X_\alpha$. Then the skeletality of the limit map $f:X\to Y$ yields an open set $V_U\subset Y$ such that $V_U\subset\cl f(p_\alpha^{-1}(U))$. By the definition of the topology of the limit space $Y$, for some index $\alpha_U\in A$, $\alpha_U\ge\alpha$, there is a non-empty open set $W_U\subset Y_{\alpha_U}$ such that $\pi_{\alpha_U}^{-1}(W_U)\subset V_U$. Since the index set $A$ is $\tau$-directed, the set $\{\alpha_U:U\in\mathcal B\}$ has an upper bound $\beta$ in $A$. It is easy to see that the index $\beta$ has the property stated in Claim~\ref{cl:1nn}.
 \end{proof}

\begin{claim}\label{cl:2nn} The set $B$ is cofinal in $A$.
\end{claim}

\begin{proof} Fix any index $\alpha_0\in A$. Using Claim~\ref{cl:1nn}, by induction we can construct a non-decreasing sequence $(\alpha_n)_{n\in\w}$ in $A$ such that
 for any non-empty open set $U\subset X_{\alpha_n}$, $n\in\w$,  there is a non-empty open set $V\subset Y_{\alpha_{n+1}}$ with  $\pi_{\alpha_{n+1}}^{-1}(V)\subseteq\cl f(p_{\alpha_n}^{-1}(U))$.

Since the set $A$ is $\w$-complete, the set $\{\alpha_n\}_{n\in\w}$ has the smallest upper bound $\beta=\sup\{\alpha_n\}_{n\in\w}\in A$. The proof of Claim~\ref{cl:2nn} will be complete as soon as we check that $\beta\in B$, which means that the limit $\downarrow_\beta$-square is skeletal.

Given any non-empty open set $U_\beta\subset X_\beta$ we need to find a non-empty open set $V_\beta\subset Y_\beta$ such that $\pi_\beta^{-1}(V_\beta)\subset\cl f(p_\beta^{-1}(U_\beta))$. Since the spectrum $\mathcal S_X$ is $\tau$-continuous, the space $X_\beta$ can be identified with the limit of the inverse spectrum $\{X_{\alpha_n},p_{\alpha_n}^{\alpha_m},\w\}$ and hence for the open set $U_\beta\subset X_\beta$ there are an index $n\in\IN$ and a non-empty open set $U\subset X_{\alpha_n}$ such that $(p^\beta_{\alpha_n})^{-1}(U)\subset U_\beta$. By the construction of the sequence $(\alpha_k)$, for the set $U\subset X_{\alpha_n}$  there is a non-empty open set $V\subset Y_{\alpha_{n+1}}$ such that
$\pi^{-1}_{\alpha_{n+1}}(V)\subset \cl f(p_{\alpha_n}^{-1}(U))$.

Consider the open set $V_\beta=(\pi^\beta_{\alpha_{n+1}})^{-1}(V)\subset Y_\beta$. Taking into account that the limit projections $p_\beta$ and $\pi_\beta$ are surjective, we conclude that
$$\begin{aligned}
V_\beta&=\pi_\beta(\pi^{-1}_{\alpha_{n+1}}(V))\subset\pi_\beta(\cl f(p_{\alpha_n}^{-1}(U))\subset
\cl \pi_\beta\circ f(p_{\alpha_n}^{-1}(U))=\\
&=\cl f_\beta\circ p_\beta(p_{\alpha_n}^{-1}(U))\subset
\cl f_\beta((p_{\alpha_n}^\beta)^{-1}(U))\subset \cl f_\beta(U_\beta),
\end{aligned}
$$
witnessing that $\beta\in B$.
\end{proof}

\begin{claim} The set $B$ is $\tau$-closed in $A$.
\end{claim}

\begin{proof} Let $C\subset B$ be a directed subset of cardinality $|C|\le\tau$ having the lowest upper bound $\gamma=\sup C$ in $A$. We need to show that $\gamma\in B$, which means that the limit $\gamma$-square is skeletal. Fix a non-empty open subset $U_\gamma\subset X_\gamma$.
Since the spectrum $\mathcal S_X$ is $\tau$-continuous, the space $X_\gamma$ can be identified with the limit space of the inverse spectrum $\{X_{\alpha},p_{\alpha}^{\beta},C\}$. Then the open set $U_\gamma\subset X_\gamma$ contains
the preimage $(p_{\alpha}^\gamma)^{-1}(U_\alpha)$ of some non-empty open set $U_\alpha\subset X_{\alpha}$, $\alpha\in C$. Since $\alpha\in C\subset B$, the limit $\downarrow_{\alpha}$-square is skeletal. Consequently, for the set $U_\alpha$ there is a non-empty open set $V_\alpha\subset Y_\alpha$ such that
$\pi_\alpha^{-1}(V_\alpha)\subset\cl f(p_\alpha^{-1}(U_\alpha))$. Then for the open subset $V_\gamma=(\pi^\gamma_\alpha)^{-1}(V_\alpha)$ in $X_\gamma$ we get
$$\pi_\gamma^{-1}(V_\gamma)=\pi_\alpha^{-1}(V_\alpha)\subset \cl f(p_\alpha^{-1}(U_\alpha))=\cl f\big(p_\gamma^{-1}((p_\alpha^\gamma)^{-1}(U_\alpha))\big)\subset \cl f(p_\gamma^{-1}(U_\gamma)),$$
witnessing that the limit $\downarrow_\gamma$-square is skeletal.
 \end{proof}

\begin{claim} For any indices $\alpha\le\beta$ in $B$ the bonding $\downarrow_\alpha^\beta$-square is skeletal.
\end{claim}

\begin{proof}  To show that the  bonding $\downarrow_\alpha^\beta$-square is skeletal, fix any open non-empty subset $U\subseteq X_\alpha$. Since $\alpha\in B$, the limit $\downarrow_\alpha$-square is skeletal and hence there exists open non-empty subset $V\subseteq Y_\alpha$ such that
$\pi_\alpha^{-1}(V)\subseteq \cl f(p_\alpha^{-1}(U))$.
Since the limit projections $p_\beta$ and $\pi_\beta$ are surjective, we get
$$(\pi^\beta_\alpha)^{-1}(V)=\pi_\beta(\pi_\alpha^{-1}(V))\subseteq \pi_\beta(\cl f(p_\alpha^{-1}(U)))
\subseteq \cl \pi_\beta\circ f(p_\alpha^{-1}(U))=\cl f_\beta\circ p_\beta(p_\alpha^{-1}(U))=\cl f_\beta((p^\beta_\alpha)^{-1}(U)).$$
\end{proof}

The definition of the set $B$ and Remark~\ref{rem1} imply our last claim, which completes the proof of Theorem~\ref{p2}.

\begin{claim}
For every $\alpha\in B$ the map $f_\alpha:X_\alpha\to Y_\alpha$ is skeletal and hence the morphism $\{f_\alpha\}_{\alpha\in B}$ is skeletal.
\end{claim}
\end{proof}

The following theorem partly reverses Theorem~\ref{p2}.

\begin{theorem}\label{p2a} Let $\{f_\alpha\}_{\alpha\in A}:\mathcal S_X\to \mathcal S_Y$ be a morphism between two $\pi\tau$-spectra $\mathcal S_X=\{X_\alpha,p^\beta_\alpha,A\}$ and $\mathcal S_Y=\{Y_\alpha,\pi^\beta_\alpha,A\}$ with surjective limit projections. If the limit map $\lim f_\alpha:\lim \mathcal S_X\to\lim\mathcal S_Y$ is not skeletal, then the set $$B=\{\alpha\in A:f_\alpha\mbox{ is not skeletal}\}$$is  $\w$-stationary in $A$.
\end{theorem}

\begin{proof} Assume that the limit map $f=\lim f_\alpha:X\to Y$ between the limit spaces $X=\lim\mathcal S_X$ and $Y=\lim\mathcal S_Y$ is not skeletal. Then the space $X$ contains a non-empty open set $U\subset V$ whose image $f(U)$ is nowhere dense in $Y$. We lose no generality assuming that the set $U$ is of the form $U=p_o^{-1}(U_o)$ for some index $o\in A$ and some non-empty open set $U_o\subset X_o$.

To prove our theorem, we need to check that the set $B$  meets each cofinal $\w$-closed subset $C$ of $A$.

\begin{claim}\label{cl:1nnn} For any index $\alpha\in C$, $\alpha\ge o$, there is an index $\beta\in C$, $\beta\ge\alpha$, such that for any non-empty open set $V_\alpha\subset Y_\alpha$ there is a non-empty open set $W_\beta\subset Y_\beta$ such that $\pi_\beta^{-1}(W_\beta)\subset \pi_\alpha^{-1}(V_\alpha)\setminus f(U)$.
\end{claim}

\begin{proof} Fix a $\pi$-base $\mathcal B$ for the space $Y_\alpha$ having cardinality $|\mathcal B|=\pi w(Y_\alpha)\le\kappa$. Since the set $f(U)$ is nowhere dense, for every set $V\in\mathcal B$ the open subset $\pi_\alpha^{-1}(V)\setminus\cl f(U)$ of $Y$ is not empty and hence contains a  set of the form $\pi_{\alpha_V}^{-1}(W_V)$ for some index $\alpha_V\ge\alpha$ in $A$ and some non-empty open set $W_V\subset Y_{\alpha_V}$. Since the index set $A$ is $\kappa$-directed and the set $C$ is cofinal in $A$, the set $\{\alpha_V:V\in\mathcal B\}$ has an upper bound $\beta\in C$. It is easy to see that the index $\beta$ has the required property.
\end{proof}

Using Claim~\ref{cl:1nnn}, by induction construct a non-decreasing sequence $(\alpha_n)_{n\in\w}$ in $C$ such that $\alpha_0\ge o$ and for any non-empty open set $V\subset Y_{\alpha_n}$, $n\in\w$, there is a non-empty open set $W\subset Y_{\alpha_{n+1}}$ such that $\pi_{\alpha_{n+1}}^{-1}(W)\subset\pi_{\alpha_n}^{-1}(V)\setminus f(U)$.

Since the set $C$ is $\w$-closed in the $\w$-complete set $A$ the chain $\{\alpha_n\}_{n\in\w}\subset C$ has an lowest upper bound $\beta\in A$, which belongs to the $\w$-closed set $C$.

\begin{claim} $\beta\in B\cap C$.
\end{claim}

\begin{proof} We need to show that the map $f_\beta:X_\beta\to Y_\beta$ is not skeletal. Assuming the opposite, for the non-empty open subset $U_\beta=(p_0^\beta)^{-1}(U_o)=p_\beta(U)$ of $X_\beta$, we can find a non-empty open set $V_\beta\subset Y_\beta$ that lies in the closure $\cl f_\beta(U_\beta)$. Since the spectrum $\mathcal S_Y$ is $\w$-continuous, the space $Y_\beta$ can be identified with the limit space of the inverse spectrum $\{Y_{\alpha_n},\pi_{\alpha_n}^{\alpha_m},\w\}$. Therefore, we lose no generality assuming that the set $V_\beta$ is of the form $V_\beta=(\pi_{\alpha_n}^\beta)^{-1}(V)$ for some open set $V\subset Y_{\alpha_n}$, $n\in\w$. By the choice of $\alpha_n$, there is a non-empty open set $W\subset Y_{\alpha_{n+1}}$ such that $\pi_{\alpha_{n+1}}^{-1}(W)\subset \pi_{\alpha_n}^{-1}(V)\setminus f(U)$. Applying to this inclusion the surjective map $\pi_\beta$, we obtain that the non-empty open subset
$$
\begin{aligned}
(\pi_{\alpha_{n+1}}^\beta)^{-1}(W)&
=\pi_\beta(\pi_{\alpha_{n+1}}^{-1}(W))\subset\pi_\beta(\pi_{\alpha_n}^{-1}(V)\setminus f(U))=\\
&=\pi_\beta(\pi_{\alpha_n}^{-1}(V))\setminus\pi_\beta\circ f(U)=(\pi_{\alpha_n}^\beta)^{-1}(V)\setminus f_\beta\circ p_\beta(U)=V_\beta\setminus f_\beta(U_\beta)
\end{aligned}
$$
of $V$ does not intersect the set $f_\beta(U_\beta)$ and hence cannot lie in its closure. This contradiction shows that the map $f_\beta$ is not skeletal and hence $\beta\in B\cap C$.
\end{proof}
\end{proof}

\section{A spectral characterization of skeletal maps between realcompact spaces}

In this section we prove Theorem~\ref{main} which characterizes skeletal maps between realcompact spaces and is the main result of this paper. This characterization has been applied in the paper \cite{BKK} detecting functors that preserve skeletal maps between compact Hausdorff spaces.

Let us recall that a Tychonoff space $X$ is called {\em realcompact} if each $C$-embedding $f:X\to Y$ into a Tychonoff space $Y$ is a closed embedding. An embedding $f:X\to Y$ is called a {\em $C$-embedding} if each continuous function $\varphi:f(X)\to\Rea$ has a continuous extension $\bar\varphi:Y\to\Rea$. By Theorem 3.11.3 \cite{En}, a topological space is realcompact of and only if it is homeomorphic to a closed subspace of some power ${\Rea}^\kappa$ of the real line, see \cite[\S3.11]{En}. By \cite[3.11.12]{En}, each Lindel\"of space is realcompact.

We say that two maps $f:X\to Y$ and $f':X'\to Y'$ are {\em homeomorphic} if there are homeomorphisms $h_X:X\to X'$ and $h_Y:Y\to Y'$ such that $f'\circ h_X=h_Y\circ f$. It is clear that a map $f:X\to Y$ is skeletal if and only if it is homeomorphic to a skeletal map $f':X'\to Y'$.

\begin{theorem}\label{main} For a map $f:X\to Y$ between Tychonoff spaces the following conditions are equivalent:
\begin{enumerate}
\item $f$ is skeletal and the spaces $X,Y$ are realcompact.
\item $f$ is homeomorphic to the limit map $\lim f_\alpha:\lim \mathcal S_X\to\lim\mathcal S_Y$ of a skeletal morphism $\{f_\alpha\}:\mathcal S_X\to\mathcal S_Y$ between two  $\w$-spectra $\mathcal S_X=\{X_\alpha,p_\alpha^\beta,A\}$ and $\mathcal S_Y=\{Y_\alpha,\pi_\alpha^\beta,A\}$ with surjective limit projections.
\item $f$ is homeomorphic to the limit map $\lim f_\alpha:\lim \mathcal S_X\to\lim\mathcal S_Y$ of a morphism $\{f_\alpha\}:\mathcal S_X\to\mathcal S_Y$ with skeletal limit squares between two  $\w$-spectra $\mathcal S_X=\{X_\alpha,p_\alpha^\beta,A\}$ and $\mathcal S_Y=\{Y_\alpha,\pi_\alpha^\beta,A\}$  with surjective limit projections.
\item $f$ is homeomorphic to the limit map $\lim f_\alpha:\lim \mathcal S_X\to\lim\mathcal S_Y$ of a morphism $\{f_\alpha\}:\mathcal S_X\to\mathcal S_Y$ with skeletal bonding squares between two  $\w$-spectra $\mathcal S_X=\{X_\alpha,p_\alpha^\beta,A\}$ and $\mathcal S_Y=\{Y_\alpha,\pi_\alpha^\beta,A\}$ with surjective limit projections.
\end{enumerate}
\end{theorem}

\begin{proof} We shall prove the implications $(1)\Ra(4)\Ra(3)\Ra(2)\Ra(1)$.

$(1)\Ra(4)$ Assume that the spaces $X,Y$ are realcompact.
Then Proposition 1.3.4 and 1.3.5 of \cite{Chi} imply that the map $f$ is homeomorphic to the limit map $\lim f_\alpha:\lim\mathcal S_X\to\lim\mathcal S_Y$ of a morphism $\{f_\alpha\}_{\alpha\in A}$ between two $\w$-spectra $\mathcal S_X=\{X_\alpha,p_\alpha^\beta,A\}$ and $\mathcal S_Y=\{Y_\alpha,\pi_\alpha^\beta,A\}$ with surjective limit projections. If the map $f$ is skeletal, then Theorem~\ref{p2} yields a cofinal $\w$-bounded subset $B\subset A$ such that the morphism $\{f_\alpha\}_{\alpha\in B}$ has skeletal bonding squares. Since the set $B$ is cofinal in $A$, $f$ is homeomorphic to the limit map $\lim f_\alpha$ induced by the morphism $\{f_\alpha\}_{\alpha\in B}$ with skeletal bonding squares between the inverse $\w$-spectra $\{X_\alpha,p_\alpha^\beta,B\}$ and $\{Y_\alpha,\pi_\alpha^\beta,B\}$.

The implications $(4)\Ra(3)$ and $(3)\Ra(2)$ follow from Corollary~\ref{c1a} and  Remark~\ref{rem1}, respectively.

The final implication $(2)\Ra(1)$ follows from Theorem~\ref{p2a} and Proposition 1.3.5 \cite{Chi} saying that a Tychonoff space is homeomorphic to the limit space of an $\w$-spectrum (with surjective limit projections) if and only if it is realcompact.
\end{proof}

Let us observe that Theorem~\ref{main} does not hold for arbitrary spectra. Just take any non-skeletal map $f:X\to Y$ between zero-dimensional (metrizable) compacta and apply the following lemma.

\begin{lemma} Each continuous map $f:X\to Y$ from a topological space $X$ to a realcompact space $Y$ of covering topological dimension $\dim(Y)=0$ is homeomorphic to the limit map $\lim f_\alpha:\lim \mathcal S_X\to \mathcal S_Y$ of a skeletal morphism
$\{f_\alpha\}_{\alpha\in A}:\mathcal S_X\to \mathcal S_Y$ between inverse spectra $\mathcal S_X=\{X_\alpha,p_\alpha^\beta,A\}$ and $\mathcal S_Y=\{Y_\alpha,\pi_\alpha^\beta,A\}$.
\end{lemma}

\begin{proof} By Lemma~6.5.4 of \cite{Chi}, the zero-dimensional realcompact space $Y$ is homeomorphic to a closed subspace of the power  $\IN^\tau$ for some cardinal $\tau$. Let $A=[\tau]^{<\w}$ be the family of finite subsets of $\tau$, partially ordered by the inclusion relation. For every $\alpha\in A$ let $Y_\alpha$ be the projection of the space $Y\subset\IN^\tau$ onto the face $\IN^\alpha$ and $\pi_\alpha:Y\to Y_\alpha$ be the corresponding projection map. For any finite sets $\alpha\subset \beta$ let $\pi_\alpha^\beta:Y_\beta\to Y_\alpha$ be the corresponding bonding projection. Then the space $Y$ can be identified with the limit $\lim\mathcal S_Y$ of the inverse spectrum $\mathcal S_Y=\{Y_\alpha,\pi_\alpha^\beta,A\}$ consisting of discrete spaces $Y_\alpha$, $\alpha\in A$.

The space $X$ can be identified with the limit of the trivial spectrum $\mathcal S_X=\{X_\alpha,p_\alpha^\beta,A\}$ consisting of spaces $X_\alpha=X$ and identity bonding maps $\pi_\alpha^\beta:X_\beta\to X_\alpha$. Then the map $f$ is homeomorphic to the limit map $\lim f_\alpha:\lim\mathcal S_X\to\mathcal S_Y$ of the skeletal morphism $\{f_\alpha\}_{\alpha\in A}:\mathcal S_X\to\mathcal S_Y$ consisting of the maps $f_\alpha=\pi_\alpha\circ f:X_\alpha=X\to Y_\alpha$, $\alpha\in A$. Here we remark that each map $f_\alpha:X_\alpha\to Y_\alpha$ is skeletal because the space $Y_\alpha$ is discrete.
\end{proof}


\begin{thebibliography}{10}

\bibitem{BKK} T.~Banakh, A.~Kucharski, M.~Martynenko, {\em On functors preserving skeletal maps and skeletally generated compacta}, preprint (arXiv:1108.4197).

\bibitem{Chi} A.~Chigogidze, {\em Inverse Spectra}, Elsevier, 1996.

\bibitem{En} R.~Engelking, {\em General Topology}, Heldermann Verlag, Berlin, 1989.

\bibitem{Miod} J.~Mioduszewski, L.~Rudolf, \textit{H-closed and extremally disconnected Hausdorff spaces}, Dissert. Math. {\bf 66} (1969), 55pp.
\end{thebibliography}
\end{document}